\documentclass[12pt,english,refpage,intoc,bibliography=totoc,index=totoc,BCOR7.5mm,captions=tableheading]{amsart}
\usepackage{lmodern}

\usepackage[T1]{fontenc}
\usepackage[latin9]{inputenc}
\usepackage{geometry}
\usepackage{color}
\usepackage{babel}
\usepackage{textcomp}
\usepackage{mathrsfs}
\usepackage{amstext}
\usepackage{amsthm}
\usepackage{amssymb}
\usepackage[unicode=true,
 bookmarks=true,bookmarksnumbered=true,bookmarksopen=false,
 breaklinks=false,pdfborder={0 0 1},backref=false,colorlinks=true]
 {hyperref}
\hypersetup{pdftitle={The LyX User's Guide},
 pdfauthor={LyX Team},
 pdfsubject={LyX},
 pdfkeywords={LyX},
 linkcolor=black, citecolor=black, urlcolor=blue, filecolor=blue, pdfpagelayout=OneColumn, pdfnewwindow=true, pdfstartview=XYZ, plainpages=false}

\makeatletter
\numberwithin{equation}{section}
\numberwithin{figure}{section}
\theoremstyle{definition}
    \ifx\thechapter\undefined
      \newtheorem{defn}{\protect\definitionname}
    \else
      \newtheorem{defn}{\protect\definitionname}[chapter]
    \fi
\theoremstyle{plain}
    \ifx\thechapter\undefined
	    \newtheorem{thm}{\protect\theoremname}
	  \else
      \newtheorem{thm}{\protect\theoremname}[chapter]
    \fi
\theoremstyle{plain}
    \ifx\thechapter\undefined
  \newtheorem{cor}{\protect\corollaryname}
\else
      \newtheorem{cor}{\protect\corollaryname}[chapter]
    \fi
\theoremstyle{remark}
    \ifx\thechapter\undefined
      \newtheorem{rem}{\protect\remarkname}
    \else
      \newtheorem{rem}{\protect\remarkname}[chapter]
    \fi
\theoremstyle{plain}
    \ifx\thechapter\undefined
      \newtheorem{prop}{\protect\propositionname}
    \else
      \newtheorem{prop}{\protect\propositionname}[chapter]
    \fi

\usepackage[figure]{hypcap}

\let\myTOC\tableofcontents
\renewcommand\tableofcontents{%
  \frontmatter
  \pdfbookmark[1]{\contentsname}{}
  \myTOC
  \mainmatter }

\usepackage{fancyhdr}
\@ifundefined{extrarowheight}
 {\usepackage{array}}{}
\makeatother

\providecommand{\corollaryname}{Corollary}
\providecommand{\definitionname}{Definition}
\providecommand{\propositionname}{Proposition}
\providecommand{\remarkname}{Remark}
\providecommand{\theoremname}{Theorem}

\begin{document}
\title{The Collatz Conjecture \& Non-Archimedean Spectral Theory - Part I.5
- How To Write The (Weak) Collatz Conjecture As A Contour Integral}
\author{by Maxwell C. Siegel}
\date{4 December, 2023}
\begin{abstract}
Let $q$ be an odd prime, and let $T_{q}:\mathbb{Z}\rightarrow\mathbb{Z}$
be the Shortened $qx+1$ map, defined by $T_{q}\left(n\right)=n/2$
if $n$ is even and $T_{q}\left(n\right)=\left(qn+1\right)/2$ if
$n$ is odd. The study of the dynamics of these maps is infamous for
its difficulty, with the characterization of the dynamics of $T_{3}$
being an alternative formulation of the famous \textbf{Collatz Conjecture}.
This series of papers presents a new paradigm for studying such arithmetic
dynamical systems by way of a neglected area of ultrametric analysis
which we have termed \textbf{$\left(p,q\right)$-adic analysis}, the
study of functions from the $p$-adics to the $q$-adics, where $p$
and $q$ are distinct primes. In this, the first-and-a-halfth paper
of the series, as a first application, we show that the \textbf{numen}
$\chi_{q}$ of $T_{q}$ can be used in conjunction with the \textbf{Correspondence
Principle} (CP) and classic complex-analytic tools of analytic number
theory to reformulate the study of periodic points of $T_{q}$ in
terms of a \textbf{contour integral} via an application of \textbf{Perron's
Formula} to a Dirichlet series generated by $\chi_{q}$ and the function
$M_{q}$ introduced in the first paper in this series, for which we
establish functional equations, which we use to derive their meromorphic
continuations to the left half-plane. The hypergeometric growth of
the series as $\textrm{Re}\left(s\right)\rightarrow-\infty$ seems
to preclude direct evaluation of the contour integrals via residues,
but asymptotic results may still be achievable.
\end{abstract}

\maketitle

\section*{Notation \& Preliminaries}

All of the notation used in this paper builds upon the notation used
in \cite{first blog paper}. We recall that, for any real number $x$,
we write $\mathbb{N}_{x}$ to denote the set of all integers $\geq x$
(thus, $\mathbb{N}_{0}=\left\{ 0,1,2,\ldots\right\} $; $\mathbb{N}_{1}=\left\{ 1,2,\ldots\right\} $).
For any real number $x$, $\left\lfloor x\right\rfloor $ denotes
the floor of $x$, the largest integer $\leq x$. Additionally, for
any $n\in\mathbb{N}_{0}$, we write $\#_{1}\left(n\right)$ to denote
the number of $1$s in the binary digits of $n$; $\lambda_{2}\left(n\right)$,
meanwhile, denotes the total number of binary digits of $n$ (note
that $\lambda_{2}\left(n\right)=\left\lceil \log_{2}\left(n+1\right)\right\rceil $).
As is traditional, we use $s=\sigma+it$ to denote the complex-variable
input of our Dirichlet series. We adopt the standard convention of
extending the binomial coefficient $\binom{n}{k}$ to complex values
by way of the Gamma function:
\begin{equation}
\binom{s}{k}\overset{\textrm{def}}{=}\frac{s!}{k!\left(s-k\right)!}\overset{\textrm{def}}{=}\frac{\Gamma\left(s+1\right)}{k!\Gamma\left(s-k+1\right)}
\end{equation}

Given a function $\lambda:\mathbb{N}_{0}\rightarrow\mathbb{C}$, we
write $\zeta_{\lambda}$ to denote the \textbf{Dirichlet series /
Dirichlet generating function }of $\lambda$:
\begin{equation}
\zeta_{\lambda}\left(s\right)\overset{\textrm{def}}{=}\sum_{n=1}^{\infty}\frac{\lambda\left(n\right)}{n^{s}}\label{eq:Definition of Zeta_f}
\end{equation}

We will use both Big-O notation and the Vinogradov notation ($\ll$),
with:
\begin{align}
f\left(n\right) & \ll g\left(n\right)\textrm{ as }n\rightarrow\infty\\
 & \Updownarrow\nonumber \\
f\left(n\right) & =O\left(g\left(n\right)\right)\textrm{ as }n\rightarrow\infty
\end{align}
We also use subscripts on Big-O notation to indicate where the limiting
variable is going; e.g.:
\begin{align}
f\left(x\right) & =O_{c}\left(g\left(x\right)\right)\\
 & \Updownarrow\nonumber \\
f\left(x\right) & =O\left(g\left(x\right)\right)\textrm{ as }x\rightarrow c
\end{align}
We assume the reader has some familiarity with the basic theory of
Dirichlet series and the Mellin transform. The papers of Flajolet
et. al. on Mellin transforms\textemdash particularly \cite{Flajolet - Mellin Transforms,Flajolet - Digital sums}\textemdash cover
all the necessary information. However, for ease of reference, we
state the most important results below.

\subsubsection*{Preliminaries from Complex Analysis}
\begin{defn}[Principal Value of an Integral]
Let $c\in\mathbb{R}$, and let $f\left(s\right)$ be a complex-valued
function so that, for every real number $T>0$, the integral: 
\begin{equation}
\int_{c-iT}^{c+iT}f\left(s\right)ds
\end{equation}
exists and is finite. (Note, we do \emph{not} require $\sup_{T>0}\left|\int_{c-iT}^{c+iT}f\left(s\right)ds\right|<\infty$).
Then, we write:
\begin{equation}
\textrm{PV}\left\{ \int_{c-i\infty}^{c+i\infty}f\left(s\right)ds\right\} \overset{\textrm{def}}{=}\lim_{T\rightarrow\infty}\int_{c-iT}^{c+iT}f\left(s\right)ds
\end{equation}
provided that the limit on the right exists. We call $\textrm{PV}\int_{c-i\infty}^{c+i\infty}f\left(s\right)ds$
the \textbf{principal value }of the expression $\int_{c-i\infty}^{c+i\infty}f\left(s\right)ds$.
\end{defn}
\begin{thm}[\textbf{Perron's Formula}\cite{Flajolet - Digital sums}]
Let $\lambda:\mathbb{N}_{1}\rightarrow\mathbb{C}$, let: $c>0$ lie
in the half-plane of absolute convergence of the Dirichlet Series
$\zeta_{\lambda}$. Then, for any $m\geq1$ and any real $x\geq2$:
\begin{equation}
\frac{1}{m!}\sum_{k=1}^{\left\lfloor x\right\rfloor -1}\lambda\left(k\right)\left(1-\frac{k}{n}\right)^{m}=\frac{1}{2\pi i}\int_{c-i\infty}^{c+i\infty}x^{s}\zeta_{\lambda}\left(s\right)\frac{ds}{s\left(s+1\right)\cdots\left(s+m\right)}\label{eq:Perron's Formula, general case}
\end{equation}
When $m=0$, we have:
\begin{equation}
\sum_{k=1}^{\left\lfloor x\right\rfloor -1}\lambda\left(k\right)=\textrm{PV}\left\{ \frac{1}{2\pi i}\int_{c-i\infty}^{c+i\infty}x^{s}\zeta_{\lambda}\left(s\right)\frac{ds}{s}\right\} \label{eq:Peron's formula; summatory case}
\end{equation}
\end{thm}
Perron's Formula is a specific case of the following more general
inversion formula:
\begin{thm}[\textbf{Mellin's Summation Formula} \cite{Flajolet - Digital sums}]
Consider a function $f:\left[0,\infty\right)\rightarrow\mathbb{C}$,
and let:
\begin{equation}
\mathscr{M}\left\{ f\right\} \left(s\right)\overset{\textrm{def}}{=}\int_{0}^{\infty}x^{s-1}f\left(x\right)dx
\end{equation}
denote the Mellin transform of $f$. Consider sequences $\left\{ \lambda_{k}\right\} _{k\geq1}$
and $\left\{ \mu_{k}\right\} _{k\geq1}$ of complex and positive-real
numbers $\lambda_{k}$ and $\mu_{k}$, respectively. Then, for any
$c\in\mathbb{R}$ so that both $\mathscr{M}\left\{ f\right\} \left(s\right)$
and:
\begin{equation}
\sum_{k=1}^{\infty}\frac{\lambda_{k}}{\mu_{k}^{s}}
\end{equation}
converge absolutely on the strip $\textrm{Re}\left(s\right)=c$, we
have:
\begin{equation}
\sum_{k=1}^{\infty}\lambda_{k}f\left(\mu_{k}x\right)=\frac{1}{2\pi i}\int_{c-i\infty}^{c+i\infty}x^{-s}\left(\sum_{k=1}^{\infty}\frac{\lambda_{k}}{\mu_{k}^{s}}\right)\mathscr{M}\left\{ f\right\} \left(s\right)ds\label{eq:Mellin's Summation Formula}
\end{equation}
\end{thm}

\section{Introduction}

In \cite{first blog paper}, we proved the \textbf{Correspondence
Principle }(CP), which\textemdash among other things\textemdash showed
that periodic points of $T_{q}$ in $\mathbb{Z}\backslash\left\{ 0\right\} $
were precisely those rational integer values which $T_{q}$'s numen,
$\chi_{q}$, attained over the set of rational $2$-adic integers
which infinitely many digits. \cite{first blog paper} gave several
versions of this result. The one we shall use in this paper is:
\begin{cor}[\textbf{Correspondence Principle (Periodic Points), Ver. 2}]
\label{cor:CP v2}

\vphantom{}

I. Let $\Omega\subseteq\mathbb{Z}$ be any cycle of $T_{q}$. Then,
viewing $\Omega$ as a subset of $\mathbb{Z}_{q}$, the intersection
$\chi_{q}\left(\mathbb{Z}_{2}\right)\cap\Omega$ is non-empty. Moreover,
for every $x\in\chi_{q}\left(\mathbb{Z}_{2}\right)\cap\Omega$, there
is an $n\in\mathbb{N}_{1}$ so that: 
\begin{equation}
x=\frac{\chi_{q}\left(n\right)}{1-M_{q}\left(n\right)}
\end{equation}

\vphantom{}

II. Let $n\in\mathbb{N}_{1}$. If the quantity $x$ given by: 
\begin{equation}
x=\frac{\chi_{q}\left(n\right)}{1-M_{q}\left(n\right)}
\end{equation}
is a $2$-adic integer, then $x$ is a periodic point of $T_{q}$
in $\mathbb{Z}_{2}$; if $x$ is in $\mathbb{Z}$, then $x$ is a
periodic point of $T_{q}$ in $\mathbb{Z}$. Moreover, if $x\in\mathbb{Z}$,
then $x$ is positive if and only if $M_{q}\left(n\right)<1$, and
$x$ is negative if and only if $M_{q}\left(n\right)>1$.
\end{cor}
\begin{rem}
When we speak of $T_{q}$ on $\mathbb{Z}_{2}$, we mean the standard
extention of $T_{q}$ from a map on $\mathbb{Z}$ to a map on $\mathbb{Z}_{2}$
by way of the rule:
\begin{equation}
T_{q}\left(\mathfrak{z}\right)\overset{\textrm{def}}{=}\begin{cases}
\frac{\mathfrak{z}}{2} & \textrm{if }\mathfrak{z}\in2\mathbb{Z}_{2}\\
\frac{q\mathfrak{z}+1}{2} & \textrm{if }\mathfrak{z}\in1+2\mathbb{Z}_{2}
\end{cases}
\end{equation}
\end{rem}
\begin{rem}
Recall the identity:
\begin{equation}
\chi_{q}\left(B_{2}\left(n\right)\right)=\frac{\chi_{q}\left(n\right)}{1-M_{q}\left(n\right)},\textrm{ }\forall n\in\mathbb{N}_{1}
\end{equation}
where:
\begin{equation}
B_{2}\left(n\right)\overset{\textrm{def}}{=}\begin{cases}
0 & \textrm{if }n=0\\
\frac{n}{1-2^{\lambda_{2}\left(n\right)}} & \textrm{if }n\geq1
\end{cases}
\end{equation}
\end{rem}
In this regard, Version 2 of the CP tells us that the behavior of
the function $n\mapsto\chi_{q}\left(B_{2}\left(n\right)\right)$ contains
everything we could ever want to know about the periodic points of
$T_{q}$ in $\mathbb{Z}$. This would suggest that we study $\chi_{q}\circ B_{2}:\mathbb{N}_{0}\rightarrow\mathbb{Q}$
directly. While this can\emph{ }be done, it is not practical. If there
is one unifying tool in the study of $\chi_{q}$, it is the central
importance of $\chi_{q}$'s functional equations:
\begin{align}
\chi_{q}\left(2\mathfrak{z}\right) & =\frac{1}{2}\chi_{q}\left(\mathfrak{z}\right)\\
\chi_{q}\left(2\mathfrak{z}+1\right) & =\frac{q\chi_{q}\left(\mathfrak{z}\right)+1}{2}
\end{align}
Almost every useful result for $\chi_{q}$ involves these identities
in one form or another, and, to that end, $\chi_{q}\circ B_{2}$ is
much more difficult to work with, simply because it does not satisfy
functional equations of this affine-linear form. Indeed, a simple
computation using the functional equations of $\chi_{q}$ and $M_{q}$
shows that $\left(\chi_{q}\circ B_{2}\right)\left(n\right)$ satisfies
the functional equation:
\begin{align}
\left(\chi_{q}\circ B_{2}\right)\left(2n\right) & =\frac{1-M_{q}\left(n\right)}{2-M_{q}\left(n\right)}\left(\chi_{q}\circ B_{2}\right)\left(n\right)\\
\left(\chi_{q}\circ B_{2}\right)\left(2n+1\right) & =\frac{q\left(1-M_{q}\left(n\right)\right)\left(\chi_{q}\circ B_{2}\right)\left(n\right)+1}{2-qM_{q}\left(n\right)}
\end{align}
Additionally, the $1-M_{q}\left(n\right)$ denominator term makes
$\left(\chi_{q}\circ B_{2}\right)\left(n\right)$ behave erratically,
blowing up whenever $n$ makes the real number $\#_{1}\left(n\right)\ln q-\lambda_{2}\left(n\right)\ln2$
small. While it is concievable that $\chi_{q}\circ B_{2}$ could be
profitably studied in its own right, doing so will likely be difficult.
Fortunately, \textbf{Corollary \ref{cor:CP v2}} furnishes a much
more manageable alternative. All we need to do is cancel the denominator:
\begin{align}
\frac{\chi_{q}\left(n\right)}{1-M_{q}\left(n\right)} & =x\nonumber \\
 & \Updownarrow\nonumber \\
\left(1-M_{q}\left(n\right)\right)x-\chi_{q}\left(n\right) & =0
\end{align}
By the CP, the periodic points of $T_{q}$ are precisely those $x\in\mathbb{Z}\backslash\left\{ 0\right\} $
for which $\left(1-M_{q}\left(n\right)\right)x-\chi_{q}\left(n\right)=0$
for some $n\geq1$.

This is where Perron's Formula comes into play.
\begin{defn}
Define the Dirichlet series:
\end{defn}
\begin{equation}
\zeta_{M_{q}}\left(s\right)\overset{\textrm{def}}{=}\sum_{n=1}^{\infty}\frac{M_{q}\left(n\right)}{n^{s}}\label{eq:Definition of Zeta_M_q}
\end{equation}
and: 
\begin{equation}
\zeta_{\chi_{q}}\left(s\right)\overset{\textrm{def}}{=}\sum_{n=1}^{\infty}\frac{\chi_{q}\left(n\right)}{n^{s}}\label{eq:Definition of Zeta_Chi_q}
\end{equation}
where $s$ is a complex variable. Finally, for any integer $x$, we
define:
\begin{equation}
F_{q}\left(s,x\right)\overset{\textrm{def}}{=}\sum_{n=1}^{\infty}\frac{\left(1-M_{q}\left(n\right)\right)x-\chi_{q}\left(n\right)}{n^{s}}\label{eq:Definition of F_q of s and x}
\end{equation}

\begin{rem}
Observe the identity:
\begin{equation}
F_{q}\left(s,x\right)=\left(\zeta\left(s\right)-\zeta_{M_{q}}\left(s\right)\right)x-\zeta_{\chi_{q}}\left(s\right)\label{eq:F_q as a zeta function}
\end{equation}
where $\zeta\left(s\right)$ is, of course, the Riemann Zeta Function.
Moreover, as we will soon show, $\textrm{Re}\left(s\right)>\sigma_{q}$
is the abscissa of absolute convergence for $\zeta_{M_{q}}\left(s\right)$,
$\zeta_{\chi_{q}}\left(s\right)$, and $F_{q}\left(s,x\right)$ (here,
$x$ is arbitrary), where:
\begin{equation}
\sigma_{q}\overset{\textrm{def}}{=}\log_{2}\left(\frac{q+1}{2}\right)\label{eq:Definition of sigma_q}
\end{equation}
\end{rem}
Using \textbf{Perron's Formula} (\ref{eq:Peron's formula; summatory case}),
letting $0<a<1$, and letting $c>\sigma_{q}$, we can express the
summatory function of $\left(1-M_{q}\left(n\right)\right)x-\chi_{q}\left(n\right)$
as:
\begin{equation}
\sum_{k=1}^{n}\left(\left(1-M_{q}\left(k\right)\right)x-\chi_{q}\left(k\right)\right)\overset{\mathbb{C}}{=}\textrm{PV}\left\{ \frac{1}{2\pi i}\int_{c-i\infty}^{c+i\infty}\left(n+a\right)^{s}F_{q}\left(s,x\right)\frac{ds}{s}\right\} ,\textrm{ }\forall n\geq2
\end{equation}
Subtracting the $\left(n-1\right)$th case from the $n$th case yields:
\begin{equation}
\left(1-M_{q}\left(n\right)\right)x-\chi_{q}\left(n\right)\overset{\mathbb{C}}{=}\textrm{PV}\left\{ \frac{1}{2\pi i}\int_{c-i\infty}^{c+i\infty}\frac{\left(n+a\right)^{s}-\left(n+a-1\right)^{s}}{s}F_{q}\left(s,x\right)ds\right\} ,\textrm{ }\forall n\geq3\label{eq:Correspondence Principal as a Contour Integral}
\end{equation}

Noting that:
\begin{align}
\chi_{q}\left(B_{2}\left(1\right)\right) & =\frac{\chi_{q}\left(1\right)}{1-M_{q}\left(1\right)}=\frac{\frac{1}{2}}{1-\frac{q}{2}}=\frac{1}{2-q}\\
\chi_{q}\left(B_{2}\left(2\right)\right) & =\frac{\chi_{q}\left(2\right)}{1-M_{q}\left(2\right)}=\frac{\frac{1}{4}}{1-\frac{q}{4}}=\frac{1}{4-q}
\end{align}
upon applying \textbf{Corollary \ref{cor:CP v2}}, we obtain the paper's
main result, a reformulation of the \textbf{Weak Collatz Conjecture}
(the assertion that $1,2,4$ are the only periodic points of the Collatz
map in the set of positive integers) as a contour integral:
\begin{thm}[\textbf{Contour Integral reformulation of the Weak Collatz Conjecture}]
\label{thm:Collatz as a Contour Integral}Let $q$ be an odd prime,
and let $c>\sigma_{q}$. Then:

\vphantom{}

I. An integer $x\in\mathbb{Z}\backslash\left\{ -1,0,1\right\} $ is
a periodic point of $T_{3}$ if and only if:
\begin{equation}
\textrm{PV}\left\{ \frac{1}{2\pi i}\int_{c-i\infty}^{c+i\infty}\frac{\left(n+a\right)^{s}-\left(n+a-1\right)^{s}}{s}F_{3}\left(s,x\right)ds\right\} =0
\end{equation}
for some integer $n\geq3$ and some real number $a\in\left(0,1\right)$.

\vphantom{}

II. An integer $x\in\mathbb{Z}\backslash\left\{ -1,0\right\} $ is
a periodic point of $T_{5}$ if and only if:
\begin{equation}
\textrm{PV}\left\{ \frac{1}{2\pi i}\int_{c-i\infty}^{c+i\infty}\frac{\left(n+a\right)^{s}-\left(n+a-1\right)^{s}}{s}F_{5}\left(s,x\right)ds\right\} =0
\end{equation}
for some integer $n\geq3$ and some real number $a\in\left(0,1\right)$.

\vphantom{}

III. If $q\geq7$, an integer $x\in\mathbb{Z}\backslash\left\{ ,0\right\} $
is a periodic point of $T_{q}$ if and only if:
\begin{equation}
\textrm{PV}\left\{ \frac{1}{2\pi i}\int_{c-i\infty}^{c+i\infty}\frac{\left(n+a\right)^{s}-\left(n+a-1\right)^{s}}{s}F_{q}\left(s,x\right)ds\right\} =0
\end{equation}
for some integer $n\geq3$ and some real number $a\in\left(0,1\right)$.
\end{thm}
This Theorem is an excellent example of the power and versatility
of the numen formalism presented in \cite{first blog paper,my dissertation}.
By having reformulated the study of $T_{q}$'s dynamics in terms of
the value-distribution of $\chi_{q}$\textemdash a fundamentally analytic
problem\textemdash we have embedded our study in classical analysis
and have thereby made it accessible to all the techniques and possibilities
that analysis provides. Moreover, because of the broad applicability
of the numen formalism to the dynamics of Hydra maps, we can achieve
similar contour-integral reformulations of the question ``is $x$
a periodic point?'' to all other Hydra maps on $\mathbb{Z}$, $\mathbb{Z}^{d}$
for $d\geq2$, and their isomorphic equivalents on the rings of integers
of number fields.

In this paper, we shall analyze the Dirichlet series defined above,
deriving functional equations for them and, from those, establishing
their continuation to meromorphic functions on the complex plane.
These functions will have a half-lattice of poles in the half-plane
$\textrm{Re}\left(s\right)\leq\sigma_{q}$.

\section{An Analysis of the Dirichlet Series of $\chi_{q}$ and $M_{q}$}

As stated in the introduction, everything we will do here will end
up coming back to the functional equations for $M_{q}$ and $\chi_{q}$.
First, let us examine the abscissa of absolute convergence. To do
that, we need to compute the summatory functions of $\chi_{q}$ and
$M_{q}$.
\begin{prop}
\label{prop:The-summatory-functions}The summatory functions of $\chi_{q}$
and $M_{q}$ satisfy:

\begin{equation}
\sum_{n=0}^{2^{N}-1}\chi_{q}\left(n\right)=\begin{cases}
\frac{N}{4}2^{N} & \textrm{if }q=3\\
2^{N}\frac{\left(\frac{q+1}{4}\right)^{N}-1}{q-3} & \textrm{if }q\geq5
\end{cases}\label{eq:Summatory function of Chi_q}
\end{equation}

\begin{equation}
\sum_{n=0}^{2^{N}-1}M_{q}\left(n\right)=\frac{q}{q-1}\left(\frac{q+1}{2}\right)^{N}-\frac{1}{q-1}\label{eq:Summatory function of M_q}
\end{equation}
\end{prop}
Proof: We use an iterative-recursive method based on the functional
equations of $\chi_{q}$ and $M_{q}$. This method is fundamental,
and will be used in varying guises throughout this series of papers.

First, define: 
\begin{equation}
S_{q}\left(N\right)\overset{\textrm{def}}{=}\sum_{n=0}^{2^{N}-1}\chi_{q}\left(n\right)
\end{equation}
We then split the index of summation modulo $2$ and apply $\chi_{q}$'s
functional equations:
\begin{align*}
S_{q}\left(N\right) & =\sum_{n=0}^{2^{N-1}-1}\left(\chi_{q}\left(2n\right)+\chi_{q}\left(2n+1\right)\right)\\
 & =\sum_{n=0}^{2^{N-1}-1}\left(\frac{\chi_{q}\left(n\right)}{2}+\frac{q\chi_{q}\left(n\right)+1}{2}\right)\\
 & =\sum_{n=0}^{2^{N-1}-1}\frac{1}{2}+\frac{q+1}{2}\underbrace{\sum_{n=0}^{2^{N-1}-1}\chi_{q}\left(n\right)}_{S_{q}\left(N-1\right)}
\end{align*}
This gives us the recursion relation:
\begin{equation}
S_{q}\left(N\right)=2^{N-2}+\frac{q+1}{2}S_{q}\left(N-1\right)\label{eq:First recursion relation}
\end{equation}
Nesting this gives:
\begin{align*}
S_{q}\left(N\right) & =2^{N-2}+\frac{q+1}{2}S_{q}\left(N-1\right)\\
\left(\textrm{Use (\ref{eq:First recursion relation}) for }S_{q}\left(N-1\right)\right); & =2^{N-2}+2^{N-3}\frac{q+1}{2}+\left(\frac{q+1}{2}\right)^{2}S_{q}\left(N-2\right)\\
\left(\textrm{Use (\ref{eq:First recursion relation}) for }S_{q}\left(N-2\right)\right); & =2^{N-2}+2^{N-3}\frac{q+1}{2}+2^{N-4}\left(\frac{q+1}{2}\right)^{2}+\left(\frac{q+1}{2}\right)^{3}S_{q}\left(N-3\right)\\
 & \vdots\\
 & =\left(\frac{q+1}{2}\right)^{N}S_{q}\left(0\right)+\sum_{k=0}^{N-1}2^{N-2-k}\left(\frac{q+1}{2}\right)^{k}
\end{align*}
Since:
\begin{equation}
S_{q}\left(0\right)=\sum_{n=0}^{2^{0}-1}\chi_{q}\left(n\right)=\chi_{q}\left(0\right)=0
\end{equation}
we conclude that:
\begin{equation}
S_{q}\left(N\right)=\sum_{k=0}^{N-1}2^{N-2-k}\left(\frac{q+1}{2}\right)^{k}=2^{N-2}\sum_{k=0}^{N-1}\left(\frac{q+1}{4}\right)^{k}=\begin{cases}
\frac{N}{4}2^{N} & \textrm{if }q=3\\
2^{N}\frac{\left(\frac{q+1}{4}\right)^{N}-1}{q-3} & \textrm{if }q\geq5
\end{cases}
\end{equation}

Now, let:
\begin{equation}
S_{q}\left(N\right)\overset{\textrm{def}}{=}\sum_{n=0}^{2^{N}-1}M_{q}\left(n\right)
\end{equation}
Then, applying the same method as above, using the functional equations:
\begin{align}
M_{q}\left(2n\right) & =\frac{M_{q}\left(n\right)}{2},\textrm{ if }n\geq1\\
M_{q}\left(2n+1\right) & =\frac{qM_{q}\left(n\right)}{2},\textrm{ if }n\geq0
\end{align}
we obtain:
\begin{align*}
S_{q}\left(N\right) & =1+\overbrace{\frac{1}{2}\sum_{n=1}^{2^{N-1}-1}M_{q}\left(n\right)}^{\textrm{evens}}+\overbrace{\frac{q}{2}\sum_{n=0}^{2^{N-1}-1}qM_{q}\left(n\right)}^{\textrm{odds}}\\
 & =1+\frac{1}{2}\left(S_{q}\left(N-1\right)-M_{q}\left(0\right)\right)+\frac{q}{2}S_{q}\left(N-1\right)\\
\left(M_{q}\left(0\right)=1\right); & =\frac{1}{2}+\frac{q+1}{2}S_{q}\left(N-1\right)
\end{align*}
Nesting yields:
\begin{equation}
S_{q}\left(N\right)=\left(\frac{q+1}{2}\right)^{N}S_{q}\left(0\right)+\frac{1}{2}\sum_{n=0}^{N-1}\left(\frac{q+1}{2}\right)^{n}
\end{equation}
Here:
\begin{equation}
S_{q}\left(0\right)=\sum_{n=0}^{2^{0}-1}M_{q}\left(n\right)=M_{q}\left(0\right)=1
\end{equation}
and so, we get:
\begin{equation}
S_{q}\left(N\right)=\left(\frac{q+1}{2}\right)^{N}+\frac{1}{2}\frac{\left(\frac{q+1}{2}\right)^{N}-1}{\frac{q+1}{2}-1}=\frac{q\left(\frac{q+1}{2}\right)^{N}-1}{q-1}
\end{equation}

Q.E.D.

\vphantom{}For clarity's sake, we record the following elementary
estimate on partial sums of Dirichlet series in terms of the Big-O
behavior of the coefficients' summatory function.
\begin{prop}
\label{prop:Dirichlet series convergence Lemma}Let $c$ be an integer
$\geq2$, let $a\left(n\right)$ be a non-negative real-valued function,
and suppose that:
\begin{equation}
\sum_{n=0}^{c^{N}-1}a\left(n\right)\ll b^{N}\textrm{ as }N\rightarrow\infty
\end{equation}
for some real $b\geq1$. Then:
\begin{equation}
\left|\sum_{n=0}^{c^{N}-1}\frac{a\left(n\right)}{\left(n+1\right)^{s}}\right|\ll\left(\frac{b}{c^{\sigma}}\right)^{N}\textrm{ as }N\rightarrow\infty
\end{equation}
where $s=\sigma+it$. Consequently, the Dirichlet series converges
absolutely for $\textrm{Re}\left(s\right)>\log_{c}b$.
\end{prop}
Proof: Using \textbf{Abel's summation formula}, we have:
\begin{align*}
\sum_{n=0}^{c^{N}-1}\frac{a\left(n\right)}{\left(n+1\right)^{s}} & =\frac{1}{c^{Ns}}\sum_{n=0}^{c^{N}-1}a\left(n\right)+s\int_{1}^{c^{N}}\frac{\sum_{n=0}^{\left\lfloor x\right\rfloor -1}a\left(n\right)}{x^{s+1}}dx\\
\left(c^{y}=x;c^{y}\ln c=dx\right); & =\frac{1}{c^{Ns}}\sum_{n=0}^{c^{N}-1}a\left(n\right)+s\ln c\int_{0}^{N}\frac{\sum_{n=0}^{\left\lfloor c^{y}\right\rfloor -1}a\left(n\right)}{c^{sy}}dy
\end{align*}
Thus, for $\textrm{Re}\left(s\right)=\sigma$:
\begin{align*}
\left|\int_{0}^{N}\frac{\sum_{n=0}^{\left\lfloor c^{y}\right\rfloor -1}a\left(n\right)}{c^{sy}}dy\right| & \leq\int_{0}^{N}\frac{\sum_{n=0}^{\left\lfloor c^{y}\right\rfloor -1}a\left(n\right)}{c^{\sigma y}}dy\\
 & =\sum_{k=0}^{N-1}\int_{k}^{k+1}\frac{\sum_{n=0}^{\left\lfloor c^{y}\right\rfloor -1}a\left(n\right)}{c^{\sigma y}}dy\\
\left(a\left(n\right)\geq0\right); & \leq\sum_{k=0}^{N-1}\int_{k}^{k+1}\frac{\sum_{n=0}^{c^{k+1}-1}a\left(n\right)}{c^{\sigma y}}dy\\
 & =\sum_{k=0}^{N-1}\left(\sum_{n=0}^{2^{k+1}-1}a\left(n\right)\right)\frac{c^{-\sigma k}-c^{-\sigma\left(k+1\right)}}{\sigma\ln c}\\
 & \ll\sum_{k=0}^{N-1}b^{k+1}\frac{c^{-\sigma k}-c^{-\sigma\left(k+1\right)}}{\sigma\ln c}\\
 & \ll\sum_{k=0}^{N-1}\left(b/c^{\sigma}\right)^{k}\\
 & \ll\left(b/c^{\sigma}\right)^{N}
\end{align*}
So:
\begin{align*}
\left|\sum_{n=0}^{c^{N}-1}\frac{a\left(n\right)}{\left(n+1\right)^{s}}\right| & \leq\frac{1}{c^{N\sigma}}\sum_{n=0}^{c^{N}-1}a\left(n\right)+\sigma\ln c\int_{0}^{N}\frac{\sum_{n=0}^{\left\lfloor c^{y}\right\rfloor -1}a\left(n\right)}{c^{ys}}dy\\
 & \ll\frac{1}{c^{N\sigma}}b^{N}+\left(b/c^{\sigma}\right)^{N}\\
 & \ll\left(b/c^{\sigma}\right)^{N}
\end{align*}

Q.E.D.
\begin{prop}
Let $q$ be an odd prime, and let $x\in\mathbb{Z}$. Then, $\zeta_{\chi_{q}}\left(s\right)$
and $\zeta_{M_{q}}\left(s\right)$, and $F_{q}\left(s,x\right)$ converge
absolutely for $\textrm{Re}\left(s\right)>\sigma_{q}$.
\end{prop}
Proof: \textbf{Proposition \ref{prop:The-summatory-functions}} shows
that the summatory functions of $\chi_{q}$ and $M_{q}$ satisfy:
\[
\max\left\{ \sum_{n=0}^{2^{N}-1}\chi_{q}\left(n\right),\sum_{n=0}^{2^{N}-1}M_{q}\left(n\right)\right\} \ll\left(\frac{q+1}{2}\right)^{N}\textrm{ as }N\rightarrow\infty
\]
Applying \textbf{Proposition \ref{prop:Dirichlet series convergence Lemma}}
yields the desired result for $\zeta_{\chi_{q}}\left(s\right)$ and
$\zeta_{M_{q}}\left(s\right)$. Finally, since $F_{q}\left(s,x\right)$
is a linear combination of $\zeta_{\chi_{q}}\left(s\right)$, $\zeta_{M_{q}}\left(s\right)$,
and $\zeta\left(s\right)$, its abscissa of convergence is $\leq$
the minimum of the abscissa of convergence of those three Dirichlet
series, which, for our choice of $q$, is $\textrm{Re}\left(s\right)>\sigma_{q}$.

Q.E.D.

\vphantom{}

Next, we establish an analytic continuation of $F_{q}$ to a meromorphic
function on $\mathbb{C}$. To do this, we use the functional equations
for $M_{q}$ and $\chi_{q}$ to establish functional equations for
their ordinary generating functions (OGFs) (defined below). Exploiting
the identity:

\begin{equation}
\sum_{n=1}^{\infty}\frac{a_{n}}{n^{s}}=\frac{1}{\Gamma\left(s\right)}\int_{0}^{\infty}y^{s-1}\left(\sum_{n=1}^{\infty}a_{n}e^{-ny}\right)dy\label{eq:Mellin transform identity}
\end{equation}
for all $\textrm{Re}\left(s\right)$ greater than the abscissa of
absolute convergence of the Dirichlet series on the leftight, we then
convert the OGFs' functional equations into functional equations for
our Dirichlet series.
\begin{defn}
Define the OGFs $g_{M_{q}}$ and $g_{\chi_{q}}$ by:
\end{defn}
\begin{align}
g_{M_{q}}\left(z\right) & \overset{\textrm{def}}{=}\sum_{n=1}^{\infty}M_{q}\left(n\right)z^{n}\\
g_{\chi_{q}}\left(z\right) & \overset{\textrm{def}}{=}\sum_{n=1}^{\infty}\chi_{q}\left(n\right)z^{n}
\end{align}

\begin{prop}
\label{prop:g_q functional equations}The functional equations:
\begin{align}
g_{M_{q}}\left(z^{2}\right) & =\frac{2}{qz+1}g_{M_{q}}\left(z\right)-\frac{qz}{qz+1}\\
g_{\chi_{q}}\left(z^{2}\right) & =\frac{2}{qz+1}g_{\chi_{q}}\left(z\right)-\frac{z}{qz+1}\frac{1}{1-z^{2}}
\end{align}
hold for all complex $z$ with $\left|z\right|<1$.
\end{prop}
Proof:

I. 
\begin{align*}
g_{M_{q}}\left(z\right) & =\sum_{n=1}^{\infty}M_{q}\left(n\right)z^{n}\\
 & =\sum_{n=0}^{\infty}M_{q}\left(2n+1\right)z^{2n+1}+\sum_{n=1}^{\infty}M_{q}\left(2n\right)z^{2n}\\
 & =\frac{qz}{2}\sum_{n=0}^{\infty}M_{q}\left(n\right)z^{2n}+\frac{1}{2}\sum_{n=1}^{\infty}M_{q}\left(n\right)z^{2n}\\
 & =\frac{qz}{2}\left(1+g_{M_{q}}\left(z^{2}\right)\right)+\frac{1}{2}g_{M_{q}}\left(z^{2}\right)\\
 & =\frac{qz}{2}+\frac{qz+1}{2}g_{M_{q}}\left(z^{2}\right)
\end{align*}
Hence:
\begin{align}
\frac{qz}{2}+\frac{qz+1}{2}g_{M_{q}}\left(z^{2}\right) & =g_{M_{q}}\left(z\right)\\
 & \Updownarrow\nonumber \\
g_{M_{q}}\left(z^{2}\right) & =\frac{2}{qz+1}\left(g_{M_{q}}\left(z\right)-\frac{qz}{2}\right)
\end{align}

II.

\begin{align*}
g_{\chi_{q}}\left(z\right) & =\sum_{n=1}^{\infty}\chi_{q}\left(n\right)z^{n}\\
 & =\sum_{n=0}^{\infty}\chi_{q}\left(2n+1\right)z^{2n+1}+\sum_{n=1}^{\infty}\chi_{q}\left(2n\right)z^{2n}\\
 & =\sum_{n=0}^{\infty}\frac{q\chi_{q}\left(n\right)+1}{2}z^{2n+1}+\frac{1}{2}\sum_{n=1}^{\infty}\chi_{q}\left(n\right)z^{2n}\\
 & =\frac{qz}{2}\sum_{n=0}^{\infty}\chi_{q}\left(n\right)z^{2n}+\frac{1}{2}\sum_{n=0}^{\infty}z^{2n+1}+\frac{1}{2}g_{\chi_{q}}\left(z^{2}\right)\\
\left(\chi_{q}\left(0\right)=0\right); & =\frac{qz}{2}g_{\chi_{q}}\left(z^{2}\right)+\frac{1}{2}\frac{z}{1-z^{2}}++\frac{1}{2}g_{\chi_{q}}\left(z^{2}\right)\\
 & =\frac{1}{2}\frac{z}{1-z^{2}}+\frac{qz+1}{2}g_{\chi_{q}}\left(z^{2}\right)
\end{align*}
Hence:
\begin{align}
\frac{1}{2}\frac{z}{1-z^{2}}+\frac{qz+1}{2}g_{\chi_{q}}\left(z^{2}\right) & =g_{\chi_{q}}\left(z\right)\\
 & \Updownarrow\nonumber \\
g_{\chi_{q}}\left(z^{2}\right) & =\frac{2}{qz+1}\left(g_{\chi_{q}}\left(z\right)-\frac{1}{2}\frac{z}{1-z^{2}}\right)
\end{align}

Q.E.D.

\vphantom{}To simplify matters, we now define an OGF for the sequence
$\left(1-M_{q}\left(n\right)\right)x-\chi_{q}\left(n\right)$.
\begin{defn}
For each $x\in\mathbb{Z}$, define the function $z\mapsto f_{q}\left(z,x\right)$
by:
\begin{equation}
f_{q}\left(z,x\right)\overset{\textrm{def}}{=}\sum_{n=1}^{\infty}\left(\left(1-M_{q}\left(n\right)\right)x-\chi_{q}\left(n\right)\right)z^{n}=\left(\frac{z}{1-z}-g_{M_{q}}\left(z\right)\right)x-g_{\chi_{q}}\left(z\right)
\end{equation}
\end{defn}
The functional equations for the $g$s give us a functional equation
for $f_{q}$:
\begin{prop}
We have:

\begin{equation}
f_{q}\left(z,x\right)=-\frac{\left(q-1\right)x}{2}\frac{z}{1-z}-\frac{1}{2}\frac{z}{1-z^{2}}+\frac{qz+1}{2}f_{q}\left(z^{2},x\right),\textrm{ }\forall\left|z\right|<1\label{eq:f_q functional equation}
\end{equation}
\end{prop}
Proof: Using \textbf{Proposition \ref{prop:g_q functional equations}},
we have:
\begin{align*}
f_{q}\left(z^{2},x\right) & =\left(\frac{z}{1-z}-g_{M_{q}}\left(z^{2}\right)\right)x-g_{\chi_{q}}\left(z^{2}\right)\\
 & =\left(\frac{z}{1-z}-\left(\frac{2}{qz+1}g_{M_{q}}\left(z\right)-\frac{qz}{qz+1}\right)\right)x-\left(\frac{2}{qz+1}g_{\chi_{q}}\left(z\right)-\frac{z}{qz+1}\frac{1}{1-z^{2}}\right)\\
 & =\left(\frac{z}{1-z}+\frac{qz}{qz+1}\right)x+\frac{z}{qz+1}\frac{1}{1-z^{2}}+\frac{2}{qz+1}\underbrace{\left(-xg_{M_{q}}\left(z\right)-g_{\chi_{q}}\left(z\right)\right)}_{\pm\frac{xz}{1-z}}\\
 & =\frac{xz}{1-z}+\frac{qxz}{qz+1}+\frac{z}{qz+1}\frac{1}{1-z^{2}}+\frac{2}{qz+1}\left(-\frac{xz}{1-z}+\underbrace{\left(\frac{z}{1-z}-g_{M_{q}}\left(z\right)\right)x-g_{\chi_{q}}\left(z\right)}_{f_{q}\left(z,x\right)}\right)\\
 & =\frac{xz}{1-z}+\frac{qxz}{qz+1}+\frac{z}{qz+1}\frac{1}{1-z^{2}}-\frac{2}{qz+1}\frac{xz}{1-z}+\frac{2}{qz+1}f_{q}\left(z,x\right)\\
 & =\frac{x}{1-z}-\frac{x}{qz+1}+\frac{z}{qz+1}\frac{1}{1-z^{2}}-\frac{x}{qz+1}\frac{2z}{1-z}+\frac{2}{qz+1}f_{q}\left(z,x\right)\\
 & =\frac{\left(q-1\right)xz\left(1+z\right)+z+2\left(1-z^{2}\right)f_{q}\left(z,x\right)}{\left(qz+1\right)\left(1-z^{2}\right)}
\end{align*}
and so:
\begin{align*}
f_{q}\left(z^{2},x\right) & =\frac{\left(q-1\right)xz\left(1+z\right)+z+2\left(1-z^{2}\right)f_{q}\left(z,x\right)}{\left(qz+1\right)\left(1-z^{2}\right)}\\
 & \Updownarrow\\
f_{q}\left(z,x\right) & =\frac{\left(qz+1\right)\left(1-z^{2}\right)f_{q}\left(z^{2},x\right)}{2\left(1-z^{2}\right)}-\frac{\left(q-1\right)xz\left(1+z\right)}{2\left(1-z^{2}\right)}-\frac{z}{2\left(1-z^{2}\right)}\\
 & =\frac{qz+1}{2}f_{q}\left(z^{2},x\right)-\frac{\left(q-1\right)x}{2}\frac{z}{1-z}-\frac{z}{2\left(1-z^{2}\right)}
\end{align*}

Q.E.D.

\vphantom{}

Using $f_{q}\left(z,x\right)$'s functional equation in conjunction
with the Mellin transform, we can obtain a functional equation for
$F_{q}\left(s,x\right)$.
\begin{thm}[\textbf{Functional Equation \& Analytic Continuation of $F_{q}$}]
\label{thm:F_q of s and x / functional equation}$F_{q}\left(s,x\right)$
satisfies the functional equation:
\begin{equation}
F_{q}\left(s,x\right)=-\frac{1}{2}\frac{2^{s}\left(\left(q-1\right)x+1\right)-1}{2^{s}-2^{\sigma_{q}}}\zeta\left(s\right)+\frac{q/2}{2^{s}-2^{\sigma_{q}}}\sum_{n=1}^{\infty}\left(-\frac{1}{2}\right)^{n}\binom{s+n-1}{n}F_{q}\left(s+n,x\right)\label{eq:Functional Equation for F_q}
\end{equation}
where:
\begin{equation}
\binom{s+n-1}{s-1}=\frac{1}{n!}\frac{\Gamma\left(s+n\right)}{\Gamma\left(s\right)}
\end{equation}
\end{thm}
Proof: Using the identity:which holds for all $s$ with $\textrm{Re}\left(s\right)>\sigma_{q}$,
we have:
\begin{equation}
F_{q}\left(s,x\right)=\frac{1}{\Gamma\left(s\right)}\int_{0}^{\infty}y^{s-1}f_{q}\left(e^{-y},x\right)dy
\end{equation}
Replacing $z$ with $e^{-y}$ in (\ref{eq:f_q functional equation})
and applying (\ref{eq:Mellin transform identity}) yields
\begin{align}
F_{q}\left(s,x\right) & =\underbrace{\frac{1}{\Gamma\left(s\right)}\int_{0}^{\infty}y^{s-1}\frac{qe^{-y}+1}{2}f_{q}\left(e^{-2y},x\right)dy}_{\textrm{I}}-\frac{\left(q-1\right)x}{2}\underbrace{\frac{1}{\Gamma\left(s\right)}\int_{0}^{\infty}\frac{y^{s-1}e^{-y}}{1-e^{-y}}dy}_{\textrm{II}}\label{eq:Big Mellin Transform}\\
 & -\frac{1}{2}\underbrace{\frac{1}{\Gamma\left(s\right)}\int_{0}^{\infty}\frac{y^{s-1}e^{-y}}{1-e^{-2y}}dy}_{\textrm{III}}\nonumber 
\end{align}

(I) becomes:
\begin{equation}
\frac{1}{2}\frac{1}{\Gamma\left(s\right)}\int_{0}^{\infty}y^{s-1}f_{q}\left(e^{-2y},x\right)dy+\frac{q}{2}\frac{1}{\Gamma\left(s\right)}\int_{0}^{\infty}y^{s-1}e^{-y}f_{q}\left(e^{-2y},x\right)dy
\end{equation}
which is:
\begin{equation}
\frac{1}{2^{s+1}}\underbrace{\frac{1}{\Gamma\left(s\right)}\int_{0}^{\infty}u^{s-1}f_{q}\left(e^{-u},x\right)du}_{F_{q}\left(s,x\right)}+\frac{q}{2^{s+1}}\frac{1}{\Gamma\left(s\right)}\int_{0}^{\infty}u^{s-1}e^{-u/2}f_{q}\left(e^{-u},x\right)du
\end{equation}
For the integral on the left, we expand $e^{-u/2}$ as a power series
and integrate term by term:
\begin{align*}
\frac{1}{\Gamma\left(s\right)}\int_{0}^{\infty}u^{s-1}e^{-u/2}f_{q}\left(e^{-u},x\right)du & =\sum_{n=0}^{\infty}\frac{\left(-1/2\right)^{n}}{n!}\frac{1}{\Gamma\left(s\right)}\int_{0}^{\infty}u^{s+n-1}f_{q}\left(e^{-u},x\right)du\\
 & =\sum_{n=0}^{\infty}\frac{\left(-1/2\right)^{n}}{n!}\frac{\Gamma\left(s+n\right)}{\Gamma\left(s\right)}\underbrace{\frac{1}{\Gamma\left(s+n\right)}\int_{0}^{\infty}u^{s+n-1}f_{q}\left(e^{-u},x\right)du}_{F_{q}\left(s+n,x\right)}\\
\left(\frac{\Gamma\left(s+n\right)}{n!\Gamma\left(s\right)}=\binom{s+n-1}{n}\right); & =\sum_{n=0}^{\infty}\left(-\frac{1}{2}\right)^{n}\binom{s+n-1}{n}F_{q}\left(s+n,x\right)
\end{align*}
So:
\begin{align*}
\frac{1}{\Gamma\left(s\right)}\int_{0}^{\infty}y^{s-1}\frac{qe^{-y}+1}{2}f_{q}\left(e^{-2y},x\right)dy & =\frac{F_{q}\left(s,x\right)}{2^{s+1}}+\frac{q}{2^{s+1}}\sum_{n=0}^{\infty}\left(-\frac{1}{2}\right)^{n}\binom{s+n-1}{n}F_{q}\left(s+n,x\right)\\
 & =\frac{q+1}{2^{s+1}}F_{q}\left(s,x\right)+\frac{q}{2^{s+1}}\sum_{n=1}^{\infty}\left(-\frac{1}{2}\right)^{n}\binom{s+n-1}{n}F_{q}\left(s+n,x\right)
\end{align*}
Meanwhile, (II) is:
\begin{equation}
\frac{1}{\Gamma\left(s\right)}\int_{0}^{\infty}\frac{y^{s-1}e^{-y}}{1-e^{-y}}dy=\sum_{n=1}^{\infty}\underbrace{\frac{1}{\Gamma\left(s\right)}\int_{0}^{\infty}y^{s-1}e^{-ns}dy}_{1/n^{s}}=\zeta\left(s\right)
\end{equation}
while (III) is: 
\begin{align*}
\frac{1}{\Gamma\left(s\right)}\int_{0}^{\infty}\frac{y^{s-1}e^{-y}}{1-e^{-2y}}dy & =\sum_{n=0}^{\infty}\underbrace{\frac{1}{\Gamma\left(s\right)}\int_{0}^{\infty}y^{s-1}e^{-\left(2n+1\right)y}dy}_{1/\left(2n+1\right)^{s}}\\
 & =\sum_{n=0}^{\infty}\frac{1}{\left(2n+1\right)^{s}}\\
 & =\sum_{n=1}^{\infty}\frac{1}{n^{s}}-\sum_{n=1}^{\infty}\frac{1}{\left(2n\right)^{s}}\\
 & =\left(1-2^{-s}\right)\zeta\left(s\right)
\end{align*}

Putting everything together, (\ref{eq:Big Mellin Transform}) becomes:
\begin{align}
F_{q}\left(s,x\right) & =\frac{q+1}{2^{s+1}}F_{q}\left(s,x\right)-\frac{\left(q-1\right)x}{2}\zeta\left(s\right)-\frac{1}{2}\left(1-2^{-s}\right)\zeta\left(s\right)\\
 & +\frac{q}{2^{s+1}}\sum_{n=1}^{\infty}\left(-\frac{1}{2}\right)^{n}\binom{s+n-1}{n}F_{q}\left(s+n,x\right)\nonumber 
\end{align}
and so:
\begin{align*}
F_{q}\left(s,x\right) & =\frac{-\frac{\left(q-1\right)x}{2}\zeta\left(s\right)-\frac{1}{2}\left(1-2^{-s}\right)\zeta\left(s\right)}{1-\frac{q+1}{2^{s+1}}}+\frac{\frac{q}{2^{s+1}}}{1-\frac{q+1}{2^{s+1}}}\sum_{n=1}^{\infty}\left(-\frac{1}{2}\right)^{n}\binom{s+n-1}{n}F_{q}\left(s+n,x\right)\\
 & =-\frac{1}{2}\frac{2^{s}\left(\left(q-1\right)x+1\right)-1}{2^{s}-\frac{q+1}{2}}\zeta\left(s\right)+\frac{q/2}{2^{s}-\frac{q+1}{2}}\sum_{n=1}^{\infty}\left(-\frac{1}{2}\right)^{n}\binom{s+n-1}{n}F_{q}\left(s+n,x\right)\\
 & =-\frac{1}{2}\frac{2^{s}\left(\left(q-1\right)x+1\right)-1}{2^{s}-2^{\sigma_{q}}}\zeta\left(s\right)+\frac{q/2}{2^{s}-2^{\sigma_{q}}}\sum_{n=1}^{\infty}\left(-\frac{1}{2}\right)^{n}\binom{s+n-1}{n}F_{q}\left(s+n,x\right)
\end{align*}

Q.E.D.

\vphantom{}Using (\ref{eq:Functional Equation for F_q}), we can
analytically continue $F_{q}$ as a function of $s$. The next few
results chronicle $F_{q}$'s singularities. First, however, a notation:
\begin{defn}
We define the function $S_{q}\left(s,x\right)$ by:
\begin{equation}
S_{q}\left(s,x\right)\overset{\textrm{def}}{=}\sum_{n=1}^{\infty}\left(-\frac{1}{2}\right)^{n}\binom{s+n-1}{n}F_{q}\left(s+n,x\right)\label{eq:Definition of S_q}
\end{equation}
Note that, for any $x$ and any $\delta>0$, $S_{q}\left(s,x\right)$
converges uniformly to a holomorphic function on the half plane $\textrm{Re}\left(s\right)\geq\sigma_{q}-1+\delta$.
\end{defn}
\begin{rem}
With this notation, (\ref{eq:Functional Equation for F_q}) can be
written as:
\begin{equation}
F_{q}\left(s,x\right)=\frac{1}{2^{s}-2^{\sigma_{q}}}\left(\frac{1-2^{s}\left(\left(q-1\right)x+1\right)}{2}\zeta\left(s\right)+\frac{q}{2}S_{q}\left(s,x\right)\right)\label{eq:F_q functional equation with S_q}
\end{equation}
\end{rem}
\begin{rem}
For $n\geq1$: 
\begin{equation}
\binom{s+n-1}{n}=\frac{1}{n!}\prod_{k=0}^{n-1}\left(s+k\right)
\end{equation}
is a degree $n$ polynomial in $s$ with positive rational coefficients
and a constant term of $0$, we have that its absolue value is bounded
by a constant multiple of $\left|s\right|^{n}$, where the constant
is the value of the polynomial at $s=1$; that is:
\[
\]
\begin{equation}
\left|\binom{s+n-1}{n}\right|\leq\left(\frac{1}{n!}\prod_{k=0}^{n-1}\left(1+k\right)\right)\left|s\right|^{n}=\frac{n!}{n!}\left|s\right|^{n}=\left|s\right|^{n}
\end{equation}
and so:

\begin{equation}
\left|S_{q}\left(s,x\right)\right|\leq\sum_{n=1}^{\infty}\left(-\frac{s}{2}\right)^{n}\left|F_{q}\left(s+n,x\right)\right|
\end{equation}
\end{rem}
We conclude by computing the residues of $F_{3}$. The general case
of $F_{q}$ is more complicated.
\begin{cor}[Singularities of $F_{3}$]

I. For $k\in\mathbb{Z}\backslash\left\{ 0\right\} $, $F_{3}\left(s,x\right)$
has a simple pole at \textbf{$s=1+\frac{2k\pi i}{\ln2}$}, the residue
of which is:
\[
\frac{3}{\ln16}S_{3}\left(1+\frac{2k\pi i}{\ln2},x\right)-\frac{4x+1}{\ln16}\zeta\left(1+\frac{2k\pi i}{\ln2}\right)
\]

II. $F_{3}\left(s,x\right)$ has a double pole at $s=1$, the residue
of which is:

\[
\frac{3S_{3}\left(1,x\right)}{\ln16}-\left(\left(\frac{1}{2}+\frac{\gamma}{\ln2}\right)x+\left(\frac{3}{8}+\frac{\gamma}{4\ln2}\right)\right)
\]

III. $F_{3}\left(s,x\right)$ has a simple pole at $s=0$, the residue
of which is:
\begin{equation}
-\frac{3\left(4x+1\right)}{16\ln2}
\end{equation}
More generally, $F_{3}\left(s,x\right)$ has a simple pole at $s=-n$
for all integers $n\geq0$. Letting $R_{3,n}$ denote the residue
of $F_{3}\left(s,x\right)$ at $s=-n$, we then have that the $R_{3,n}$s
satisfy the recurrence relation:
\begin{equation}
R_{3,n+1}=\frac{\frac{3\left(4x+1\right)}{16\ln2}}{1-2^{n}}\frac{1}{n}+\frac{3/2}{1-2^{n}}\sum_{k=0}^{n}2^{k}\binom{n+1}{k}R_{3,k},\textrm{ }\forall n\geq0\label{eq:R_3,n recurrence relation}
\end{equation}
\end{cor}
Proof: For convenience, here is (\ref{eq:F_q functional equation with S_q})
again:

\[
F_{q}\left(s,x\right)=\frac{1}{2^{s}-2^{\sigma_{q}}}\left(\frac{1-2^{s}\left(\left(q-1\right)x+1\right)}{2}\zeta\left(s\right)+\frac{q}{2}S_{q}\left(s,x\right)\right)
\]
Fix $x\in\mathbb{Z}$. Note that $s=\sigma_{q}$ is the largest positive
real number at which $F_{q}\left(s,x\right)$ has a singularity, and
that $S_{q}\left(s,x\right)$ is holomorphic at this value of $s$.
Using (\ref{eq:F_q functional equation with S_q}), we see that the
singularity of $F_{q}\left(s,x\right)$ at $s=\sigma_{q}$ must come
from either:
\begin{eqnarray*}
 & \frac{1-2^{s}\left(\left(q-1\right)x+1\right)}{2}\zeta\left(s\right)\\
 & \textrm{or}\\
 & \frac{1}{2^{s}-2^{\sigma_{q}}}
\end{eqnarray*}
Since $\zeta\left(s\right)$ has a simple pole at $s=1$, the limit:
\[
\lim_{s\rightarrow1}\frac{1-2^{s}\left(\left(q-1\right)x+1\right)}{2}=-\left(q-1\right)x-\frac{1}{2}
\]
and the fact that $x\in\mathbb{Z}$ tells us the only singularity
of $\frac{1-2^{s}\left(\left(q-1\right)x+1\right)}{2}\zeta\left(s\right)$
is a simple pole at $s=1$ of residue $-\left(q-1\right)x-\frac{1}{2}$.
Meanwhile, $\left(2^{s}-2^{\sigma_{q}}\right)^{-1}$ has simple poles
at $s=\sigma_{q}+\frac{2k\pi i}{\ln2}$ for all $k\in\mathbb{Z}$.

So, suppose $q=3$, so that $\sigma_{q}=1$. Then (\ref{eq:F_q functional equation with S_q})
becomes:
\begin{equation}
F_{3}\left(s,x\right)=\left(\frac{1}{2}-2^{s-1}\left(2x+1\right)\right)\frac{\zeta\left(s\right)}{2^{s}-2}+\frac{3}{2}\frac{S_{3}\left(s,x\right)}{2^{s}-2}\label{eq:F_3 functional equation, simplified}
\end{equation}
where $S_{3}\left(s,x\right)$ is holomorphic for $\textrm{Re}\left(s\right)>\sigma_{q}-1=0$.

\textbullet{} Here:
\[
\left(\frac{1}{2}-2^{s-1}\left(2x+1\right)\right)\frac{\zeta\left(s\right)}{2^{s}-2}
\]
has a double pole at $s=1$, while $\frac{3}{2}\frac{S_{3}\left(s,x\right)}{2^{s}-2}$
has a simple pole there; in total, $F_{3}\left(s,x\right)$ has a
\textbf{double pole at $s=1$}.

To compute the residue of $F_{3}\left(s,x\right)$ at $s=1$, we use
the series expansions:
\[
\zeta\left(s\right)=\frac{1}{s-1}+\gamma+O\left(\left(s-1\right)\right)\textrm{ as }s\rightarrow1
\]

\[
\frac{s-1}{2^{s}-2}=\frac{1}{\ln4}-\frac{s-1}{4}+O\left(\left(s-1\right)^{2}\right)\textrm{ as }s\rightarrow1
\]
which gives:
\begin{align*}
\left(s-1\right)^{2}F_{3}\left(s,x\right) & =\frac{1-2^{s}\left(2x+1\right)}{\ln16}+\left(\left(\frac{\gamma}{\ln16}-\frac{1}{8}\right)\left(1-2^{s}\left(2x+1\right)\right)+\frac{3S_{3}\left(s,x\right)}{\ln16}\right)\left(s-1\right)\\
 & +O\left(\left(s-1\right)^{2}\right)
\end{align*}
Since: 
\[
\textrm{Res}\left[F_{3}\left(s,x\right):1\right]=\lim_{s\rightarrow1}\frac{d}{ds}\left\{ \left(s-1\right)^{2}F_{3}\left(s,x\right)\right\} 
\]
differentiating the previous equation with respect to $s$ and taking
the limit as $s\rightarrow1$ yields the residue:
\[
\]
\[
\textrm{Res}\left[F_{3}\left(s,x\right);1\right]=\frac{3S_{3}\left(1,x\right)}{\ln16}-\left(\left(\frac{1}{2}+\frac{\gamma}{\ln2}\right)x+\left(\frac{3}{8}+\frac{\gamma}{4\ln2}\right)\right)
\]

\textbullet{} Since $S_{3}\left(s,x\right)$ and $\zeta\left(s\right)$
are holomorphic at:

\[
s=\sigma_{3}+\frac{2k\pi i}{\ln2}=1+\frac{2k\pi i}{\ln2},\textrm{ }\forall k\in\mathbb{Z}\backslash\left\{ 0\right\} 
\]
and since:
\[
\lim_{s\rightarrow1+\frac{2k\pi i}{\ln2}}\left(\frac{1}{2}-2^{s-1}\left(2x+1\right)\right)=-2x-\frac{1}{2}
\]
which is never zero (since $x\in\mathbb{Z}$) we have that the $\left(2^{s}-2\right)^{-1}$
term is the only singular part of $F_{3}\left(s,x\right)$ at $s=1+\frac{2k\pi i}{\ln2}$
for non-zero $k$, and so, $F_{3}\left(s,x\right)$ has \textbf{simple
poles at $s=1+\frac{2k\pi i}{\ln2}$} \textbf{for} $k\in\mathbb{Z}\backslash\left\{ 0\right\} $.
The residues of these poles are:
\begin{align*}
\textrm{Res}\left[F_{3}\left(s,x\right):1+\frac{2k\pi i}{\ln2}\right] & =\lim_{s\rightarrow1+\frac{2k\pi i}{\ln2}}\left(s-1-\frac{2k\pi i}{\ln2}\right)\left(\left(\frac{1}{2}-2^{s-1}\left(2x+1\right)\right)\frac{\zeta\left(s\right)}{2^{s}-2}+\frac{3}{2}\frac{S_{3}\left(s,x\right)}{2^{s}-2}\right)\\
 & =\lim_{s\rightarrow0}\underbrace{\frac{s}{2^{s}-1}}_{1/\ln2}\frac{\left(1-2^{s}\left(4x+2\right)\right)\zeta\left(s+1+\frac{2k\pi i}{\ln2}\right)+3S_{3}\left(s+1+\frac{2k\pi i}{\ln2},x\right)}{4}\\
 & =\frac{3}{\ln16}S_{3}\left(1+\frac{2k\pi i}{\ln2},x\right)-\frac{4x+1}{\ln16}\zeta\left(1+\frac{2k\pi i}{\ln2}\right)
\end{align*}

\textbullet{} As $s\rightarrow0$, (\ref{eq:F_q functional equation with S_q})
tends to:
\[
F_{3}\left(0,x\right)=\left(\frac{1}{2}-\frac{2x+1}{2}\right)\frac{\zeta\left(0\right)}{1-2}+\frac{3}{2}\frac{S_{3}\left(0,x\right)}{1-2}=-\frac{x}{2}-\frac{3}{2}S_{3}\left(0,x\right)
\]
Now:
\[
S_{3}\left(s,x\right)=-\frac{s}{2}F_{3}\left(s+1,x\right)+\frac{s\left(s+1\right)}{8}F_{3}\left(s+2,x\right)-\frac{s\left(s+1\right)\left(s+2\right)}{2^{3}\cdot3!}F_{3}\left(s+3,x\right)+\cdots
\]
Since $F_{3}\left(s+n,x\right)$ is holomorphic at $s=0$ uniformly
with respect to $n\geq2$, we have:
\[
S_{3}\left(s,x\right)=-\frac{1}{2}sF_{3}\left(s+1,x\right)+O_{0}\left(s\right)\textrm{ as }s\rightarrow0
\]
Since $F_{3}\left(s,x\right)$ has a double pole at $s=1$, $sF_{3}\left(s+1,x\right)$
has a simple pole at $s=0$, and thus, $S_{3}\left(0,x\right)$ is
a simple pole. Consequently, $F_{3}\left(s,x\right)$ has a \textbf{simple
pole at $s=0$}, the residue of which is given by:
\begin{align*}
\lim_{s\rightarrow0}sF_{3}\left(s,x\right) & =\lim_{s\rightarrow0}s\left(-\frac{x}{2}-\frac{3}{2}S_{3}\left(s,x\right)\right)\\
 & =-\frac{3}{2}\lim_{s\rightarrow0}sS_{3}\left(s,x\right)\\
\left(S_{3}\left(s,x\right)=-\frac{1}{2}sF_{3}\left(s+1,x\right)+O_{0}\left(s\right)\right); & =\frac{3}{4}\lim_{s\rightarrow0}s^{2}F_{3}\left(s+1,x\right)\\
 & =\frac{3}{4}\lim_{s\rightarrow1}\left(s-1\right)^{2}F_{3}\left(s,x\right)\\
\left(\textrm{Use }(\ref{eq:F_3 functional equation, simplified})\right); & =\frac{3}{4}\lim_{s\rightarrow1}\left(s-1\right)^{2}\left(\frac{1}{2}-2^{s-1}\left(2x+1\right)\right)\frac{\zeta\left(s\right)}{2^{s}-2}\\
 & +\frac{3}{4}\lim_{s\rightarrow1}\left(s-1\right)^{2}\frac{3}{2}\frac{S_{3}\left(s,x\right)}{2^{s}-2}
\end{align*}
Here, $S_{3}\left(s,x\right)$ is holomorphic at $s=1$. Since $1/\left(2^{s}-2\right)$
has a simple pole at $s=1$, we have:
\[
\lim_{s\rightarrow1}\left(s-1\right)^{2}\frac{3}{2}\frac{S_{3}\left(s,x\right)}{2^{s}-2}=0
\]
On the other hand, $\zeta\left(s\right)/\left(2^{s}-2\right)$ has
a double pole at $s=1$, so multiplying it by $\left(s-1\right)^{2}$
and letting $s\rightarrow1$ will give us a constant. Consequently:
\begin{align*}
\textrm{Res}\left[F_{3}\left(s,x\right):0\right] & =\lim_{s\rightarrow0}sF_{3}\left(s,x\right)\\
 & =-\frac{3}{2}\lim_{s\rightarrow0}sS_{3}\left(s,x\right)\\
 & =\frac{3}{4}\lim_{s\rightarrow1}\left(s-1\right)^{2}F_{3}\left(s,x\right)\\
 & =\frac{3}{4}\lim_{s\rightarrow1}\left(s-1\right)^{2}\overbrace{\left(\frac{1}{2}-2^{s-1}\left(2x+1\right)\right)}^{\rightarrow-2x-\frac{1}{2}}\frac{\zeta\left(s\right)}{2^{s}-2}\\
 & =-\frac{3}{4}\left(-2x-\frac{1}{2}\right)\times\lim_{s\rightarrow1}\frac{\left(s-1\right)^{2}\zeta\left(s\right)}{2^{s}-2}\\
 & =-\frac{3}{8}\left(4x+1\right)\times\lim_{s\rightarrow1}\left(s-1\right)\zeta\left(s\right)\times\lim_{z\rightarrow1}\frac{z-1}{2^{z}-2}\\
 & =-\frac{3\left(4x+1\right)}{16\ln2}
\end{align*}
Note that this also shows that:
\begin{equation}
\lim_{s\rightarrow1}\left(s-1\right)^{2}F_{3}\left(s,x\right)=-\frac{4x+1}{4\ln2}\label{eq:simple residue of F_3 at 1}
\end{equation}

More generally, by induction, (\ref{eq:F_q functional equation with S_q})
shows that $F_{3}\left(s,x\right)$ has a \textbf{simple pole at $s=-k$
for $k\geq0$}, the residue of which is:
\begin{align*}
\lim_{s\rightarrow-k}\left(s+k\right)F_{3}\left(s,x\right) & =\lim_{s\rightarrow-k}\left(s+k\right)\frac{\left(\frac{1}{2}-2^{s-1}\left(2x+1\right)\right)\zeta\left(s\right)+\frac{3}{2}S_{3}\left(s,x\right)}{2^{s}-2}\\
 & =\frac{3/2}{2^{-k}-2}\lim_{s\rightarrow-k}\left(s+k\right)S_{3}\left(s,x\right)\\
 & =\frac{3/2}{2^{-k}-2}\lim_{s\rightarrow0}sS_{3}\left(s-k,x\right)
\end{align*}
Here:
\begin{align*}
S_{3}\left(s-k,x\right) & =\sum_{n=1}^{k}\left(-\frac{1}{2}\right)^{n}\binom{s+n-k-1}{n}F_{3}\left(s+n-k,x\right)\\
 & +\left(-\frac{1}{2}\right)^{k}\sum_{n=1}^{\infty}\left(-\frac{1}{2}\right)^{n}\binom{s+n-1}{n+k}F_{3}\left(s+n,x\right)\\
 & =\sum_{n=0}^{k-1}\left(-\frac{1}{2}\right)^{k-n}\binom{s-n-1}{k-n}F_{3}\left(s-n,x\right)\\
 & +\left(-\frac{1}{2}\right)^{k}\sum_{n=1}^{\infty}\left(-\frac{1}{2}\right)^{n}\binom{s+n-1}{n+k}F_{3}\left(s+n,x\right)
\end{align*}
As $s\rightarrow0$, note that $F_{3}\left(s-n,x\right)$ tends to
the simple pole $F_{3}\left(-n,x\right)$ for all $n\in\left\{ 0,\ldots,k-1\right\} $.
Since:
\begin{align*}
\binom{s-n-1}{k-n} & =\frac{1}{\left(k-n\right)!}\prod_{j=0}^{k-n-1}\left(s-n-1-j\right)\\
\left(\ell=j+n+1\right); & =\frac{1}{\left(k-n\right)!}\prod_{\ell=n+1}^{k}\left(s-\ell\right)\\
\left(\textrm{let }s\rightarrow0\right); & =\frac{1}{\left(k-n\right)!}\prod_{\ell=n+1}^{k}\left(-\ell\right)\\
 & =\frac{\left(-1\right)^{k-n}}{\left(k-n\right)!}\frac{\prod_{\ell=1}^{k}\ell}{\prod_{\ell=1}^{n}\ell}\\
 & =\left(-1\right)^{k-n}\frac{k!}{n!\left(k-n\right)!}\\
 & =\left(-1\right)^{k-n}\binom{k}{n}
\end{align*}
we have that:
\begin{align*}
\lim_{s\rightarrow0}s\sum_{n=0}^{k-1}\left(-\frac{1}{2}\right)^{k-n}\binom{s-n-1}{k-n}F_{3}\left(s-n,x\right) & =\sum_{n=0}^{k-1}\left(-\frac{1}{2}\right)^{k-n}\left(-1\right)^{k-n}\binom{k}{n}\textrm{Res}\left[F_{3}\left(s,x\right):-n\right]\\
 & =\sum_{n=0}^{k-1}\left(\frac{1}{2}\right)^{k-n}\binom{k}{n}\textrm{Res}\left[F_{3}\left(s,x\right):-n\right]
\end{align*}
Meanwhile, since $F_{3}\left(s+n,x\right)$ is holomorphic at $s=0$
uniformly with respect to $n\geq2$, we have:
\begin{align*}
\lim_{s\rightarrow0}s\left(-\frac{1}{2}\right)^{k}\sum_{n=1}^{\infty}\left(-\frac{1}{2}\right)^{n}\binom{s+n-1}{n+k}F_{3}\left(s+n,x\right) & =\left(-\frac{1}{2}\right)^{k+1}\lim_{s\rightarrow0}s\binom{s}{k+1}F_{3}\left(s+1,x\right)\\
 & =\frac{\left(-\frac{1}{2}\right)^{k+1}\prod_{j=1}^{k}\left(-j\right)}{\left(k+1\right)!}\lim_{s\rightarrow0}s^{2}F_{3}\left(s+1,x\right)\\
 & =-\frac{\left(\frac{1}{2}\right)^{k+1}}{k+1}\lim_{s\rightarrow1}\left(s-1\right)^{2}F_{3}\left(s,x\right)\\
\left(\ref{eq:simple residue of F_3 at 1}\right); & =-\frac{\left(\frac{1}{2}\right)^{k+1}}{k+1}\times-\frac{4x+1}{4\ln2}\\
 & =\frac{4x+1}{4\ln2}\frac{1}{k+1}\left(\frac{1}{2}\right)^{k+1}
\end{align*}
Thus:
\[
\lim_{s\rightarrow0}sS_{3}\left(s-k,x\right)=\frac{4x+1}{4\ln2}\frac{1}{k+1}\left(\frac{1}{2}\right)^{k+1}+\sum_{n=0}^{k-1}\left(\frac{1}{2}\right)^{k-n}\binom{k}{n}\textrm{Res}\left[F_{3}\left(s,x\right):-n\right]
\]
and so:
\begin{align*}
\textrm{Res}\left[F_{3}\left(s,x\right):-k\right] & =\lim_{s\rightarrow-k}\left(s+k\right)F_{3}\left(s,x\right)\\
 & =\frac{3/2}{2^{-k}-2}\lim_{s\rightarrow0}sS_{3}\left(s-k,x\right)\\
 & =\frac{3/2}{2^{-k}-2}\left(\frac{4x+1}{4\ln2}\frac{1}{k+1}\left(\frac{1}{2}\right)^{k+1}+\sum_{n=0}^{k-1}\left(\frac{1}{2}\right)^{k-n}\binom{k}{n}\textrm{Res}\left[F_{3}\left(s,x\right):-n\right]\right)\\
 & =\frac{\frac{3\left(4x+1\right)}{16\ln2}}{1-2^{k+1}}\frac{1}{k+1}+\frac{3/2}{1-2^{k+1}}\sum_{n=0}^{k-1}2^{n}\binom{k}{n}\textrm{Res}\left[F_{3}\left(s,x\right):-n\right]
\end{align*}
Thus, letting:
\[
R_{3,n}\overset{\textrm{def}}{=}\textrm{Res}\left[F_{3}\left(s,x\right):-n\right]
\]
we have:
\[
R_{3,k}=\frac{\frac{3\left(4x+1\right)}{16\ln2}}{1-2^{k+1}}\frac{1}{k+1}+\frac{3/2}{1-2^{k+1}}\sum_{n=0}^{k-1}2^{n}\binom{k}{n}R_{3,n}
\]
re-indexing (swapping $k$ and $n$) gives:
\[
R_{3,n}=\frac{\frac{3\left(4x+1\right)}{16\ln2}}{1-2^{n+1}}\frac{1}{n+1}+\frac{3/2}{1-2^{n+1}}\sum_{k=0}^{n-1}2^{k}\binom{n}{k}R_{3,k}
\]
which is the desired recursive formula.

Q.E.D.

\vphantom{}

In conjunction with Perron's Formula, these residues can be used to
asymptotically analyze equation (\ref{eq:Correspondence Principal as a Contour Integral}),
however, the functional equation (\ref{eq:Functional Equation for F_q})
can be used to show that $F_{q}\left(s,x\right)$ grows hyperexponentially
as $\textrm{Re}\left(s\right)\rightarrow-\infty$, which makes it
impossible to evaluate the integral in (\ref{eq:Correspondence Principal as a Contour Integral})
exactly by shifting the contour of integration to left infinity. Nevertheless,
it is conceivable that clever explotation of the properties of Perron-type
formulae and the growth of the Riemann Zeta Function and to get something
interesting out of this. Finally, since the set-up used in this paper
can be applied to any Hydra map for which the methods of \cite{first blog paper,my dissertation},
it may be of interest to explore if there are any instances\textemdash particularly
degenerate ones\textemdash where the Mellin transforms methods of
this paper might bear useful fruit.

\section*{Acknowledgements}

This series of paper are a distillation of the highlights of the author's
PhD (Mathematics) dissertation \cite{my dissertation} done at the
University of Southern California under the obliging supervision of
Professors Sheldon Kamienny and Nicolai Haydn. Thanks must also be
given to Jeffery Lagarias, Steven J. Miller, Alex Kontorovich, Andrei
Khrennikov, K.R. Matthews, Susan Montgomery, Amy Young and all the
helpful staff of the USC Mathematics Department, and all the kindly
strangers became and acquainted with along the way.

\end{document}